\documentclass{article}

\usepackage{algorithm, algpseudocode}
\usepackage{amsfonts, amsmath, amsopn, amsthm, epsfig}
  \numberwithin{equation}{section}
\usepackage{graphicx, epstopdf}
\usepackage{hyperref}
\usepackage{subfigure}

\usepackage{calc}
\usepackage[top=2.5cm,bottom=2.3cm,inner=2.3cm,outer=2.3cm,bindingoffset=0.5cm,footskip=0.55cm]{geometry}
\newtheorem{theorem}{Theorem}[section]
\newtheorem{corollary}[theorem]{Corollary}
\newtheorem{conjecture}[theorem]{Conjecture}

\theoremstyle{remark}
\newtheorem{remark}{Remark}

\newenvironment{lemma*}[2][Lemma]{\par\bgroup{\bfseries #1\ #2. }\it\ignorespaces}{\egroup}

\title{Online and stochastic Douglas-Rachford splitting method for large scale machine learning}

\author{Ziqiang Shi\footnotemark[1] \footnotemark[2],  Rujie Liu\footnotemark[1]
}

\newcommand{\BALD}{\begin{aligned}}
\newcommand{\EALD}{\end{aligned}}
\newcommand{\BALDS}{\begin{aligned*}}
\newcommand{\EALDS}{\end{aligned*}}
\newcommand{\BCAS}{\begin{cases}}
\newcommand{\ECAS}{\end{cases}}
\newcommand{\BEAS}{\begin{eqnarray*}}
\newcommand{\EEAS}{\end{eqnarray*}}
\newcommand{\BEQ}{\begin{equation}}
\newcommand{\EEQ}{\end{equation}}
\newcommand{\BIT}{\begin{itemize}}
\newcommand{\EIT}{\end{itemize}}
\newcommand{\BMAT}{\begin{bmatrix}}
\newcommand{\EMAT}{\end{bmatrix}}
\newcommand{\BNUM}{\begin{enumerate}}
\newcommand{\ENUM}{\end{enumerate}}

\newcommand{\BA}{\begin{array}}
\newcommand{\EA}{\end{array}}

\newcommand{\reals}{\mathbf{R}}

\DeclareMathOperator*{\argmin}{\arg\min}

\DeclareMathOperator*{\minimize}{minimize}

\DeclareMathOperator{\sign}{sign}

\newcommand{\abs}[1]{\left| #1 \right|}

\newcommand{\norm}[1]{\left\| #1 \right\|}

\DeclareMathOperator{\dom}{dom}

\DeclareMathOperator{\prox}{prox}

\begin{document}

\maketitle

\renewcommand{\thefootnote}{\fnsymbol{footnote}}

\footnotetext[1]{Fujitsu Research \& Development Center, Beijing, China.}
\footnotetext[2]{shiziqiang@cn.fujitsu.com}

\renewcommand{\thefootnote}{\arabic{footnote}}

\begin{abstract}
Online and stochastic learning has emerged as powerful tool in large scale optimization. In this work, we generalize the Douglas-Rachford splitting (DRs) method for minimizing composite functions to online and stochastic settings (to our best knowledge this is the first time DRs been generalized to sequential version).
We first establish an $O(1/\sqrt{T})$ regret bound for batch DRs method. Then we proved that the online DRs splitting method enjoy an $O(1)$ regret bound and stochastic DRs splitting has a convergence rate of $O(1/\sqrt{T})$. The proof is simple and intuitive, and the results and technique can be served as a initiate for the research on the large scale machine learning employ the DRs method. Numerical experiments of the proposed method demonstrate the effectiveness of the online and stochastic update rule, and further confirm our regret and convergence analysis.
\end{abstract}

\section{Introduction and problem statement}
\label{sec:introduction}
First introduced in~\cite{douglas1956numerical}, the Douglas-Rachford splitting technique has become popular in recent years due to its fast theoretical convergence rate and strong practical performance.
The method is first proposed to addresses
the minimization of the sum of two functions $g(x)+h(x)$.
It was extended in~\cite{lions1979splitting} to handle problems involving the sum of two nonlinear monotone operator problems. For further developments, see~\cite{eckstein1992douglas,combettes2009iterative,combettes2011proximal}. However, most of these variants implicitly assume full accessibility of all data values, while in reality one can hardly ignore the fact that the size of data is rapidly increasing in
various domain, and thus batch mode learning procedure cannot deal with the huge
size training set for the data probably cannot be loaded into the
memory simultaneously. Furthermore, it cannot be started until
the training data is prepared, hence cannot effectively deal with
training data appear in sequence, such as audio and video processing~\cite{shi2013audio}.
In such situation, sequential learning become powerful tools.

Online and stochastic learning are of the most promising methods in large scale machine learning tasks in these days~\cite{zinkevich2003online,wang2012online}.
Important advances have been made on sequential learning in the recent literature on similar problems. Composite objective mirror
descent (COMID)~\cite{duchi2010composite} generalizes mirror
descent~\cite{beck2003mirror} to the online setting.
Regularized dual averaging (RDA)~\cite{xiao2010dual} generalizes dual averaging~\cite{nesterov2009primal} to online and composite
optimization, and can be used for distributed optimization~\cite{duchi2012dual}. Online alternating direction multiplier method (ADMM)~\cite{suzuki2013dual}, RDA-ADMM~\cite{suzuki2013dual} and online proximal gradient (OPG) ADMM~\cite{wang2012online} generalize
classical ADMM~\cite{gabay1976dual} to online and stochastic settings.

Our focus in this paper is to generalize the Douglas-Rachford splitting to online and stochastic settings.
In this work, we consider the problems of the following form:
\begin{align}
  \minimize_{x \in \reals^n} \,f_T(x) := \frac{1}{T}\sum_{t=1}^T g_t (x) + h(x)=\frac{1}{T}\sum_{t=1}^T (g_t (x) + h(x)),
  \label{eq:composite-form-online}
\end{align}
where $g_t$ is a convex loss function associated with a sample in a training set, and
$h$ is a non-smooth convex penalty function or
regularizer. Many problems of relevance in signal processing and machine learning can be formulated as the above optimization problem. Similar problems include the ridge regression, the lasso~\cite{tibshirani1996regression}, the logistic regression, and the minimization of total variation.

Let $g(x)=\frac{1}{T}\sum_{t=1}^T g_t (x)$ in Problem (\ref{eq:composite-form-online}), and then the Douglas-Rachford splitting algorithm approximates a minimizer of~(\ref{eq:composite-form-online})
with the help of the following sequence:
\begin{align*}
   u_{t+1}   = u_t+\lambda_t\{\prox_{\lambda g}\left[2\prox_{\lambda h}(u_t)-u_t\right]-\prox_{\lambda h}(u_t)\},
\end{align*}
where ($\lambda_t)_{t\geq0}\subset [0, 2]$ satisfies $\sum_{t\geq1} \lambda_t(2-\lambda_t)=\infty$, and the \emph{proximal} mapping of a convex function $h$ at $x$ is
\BEQ
  \prox_{h}(t) := \argmin_{y\in\reals^n}\,h(y) + \frac{1}{2}\norm{y-t}^2. \nonumber
\EEQ
Thus the iterative scheme of Douglas-Rachford splitting for the problem~(\ref{eq:composite-form-online}) is as follows:
\begin{align} \label{eq:DRs_iterative_scheme}
   x_{t+1} % &= \prox_{\lambda_th}\left(x_t-\lambda_t\nabla g(x_t)\right)
        &= \argmin_x\,h(x) + \frac{1}{2\lambda}\norm{x-u_t}^2,
        \\ z_{t+1} & = \argmin_z\,g(z) + \frac{1}{2\lambda}\norm{z-(2x_{t+1}-u_t)}^2,
        \label{eq:DRs_iterative_scheme_z_t}
 \\   u_{t+1}  & = u_t + \lambda_t(z_{t+1}-x_{t+1}),
             \label{eq:DRs_iterative_scheme_x_t}
\end{align}
where $(u_t)_{t\in N}$ converges weakly (in Hilbert space $\reals^n$, weak convergence is equivalent to strong convergence) to some point $u$ (thus also $x_t$, for $\prox_{\lambda h}(\cdot)$ is continuous), and $\prox_{\lambda h}(u)$ is a solution to Problem (\ref{eq:composite-form-online}).
For convenience of description, we assume all $\lambda_t=1$ in this work.

The only modification of the splitting that we propose for online and stochastic processing is simple:
\begin{equation*}
        z_{t+1} = \argmin_z\,g_t(z) + \frac{1}{2\lambda}\norm{z-(2x_{t+1}-u_t)}^2,
\end{equation*}
and
\begin{equation*}
        z_{t+1} = \argmin_z\,g_{i_t}(z) + \frac{1}{2\lambda}\norm{z-(2x_{t+1}-u_t)}^2,
\end{equation*}
where index $i_t$ is sampled uniformly from the set $\{1, ..., T\}$
We call these methods online DRs (oDRs) and stochastic DRs (sDRs) respectively. Due to the complex loss function $g_t(z)$, generally the update is difficult to solve efficiently. A common way is to linearize the objective such that
\begin{equation}
        z_{t+1} = \argmin_z\,\nabla g_t(z_t)^T (z-z_t) + \frac{1}{2\lambda}\norm{z-(2x_{t+1}-u_t)}^2,
        \label{eq:inexact_Online_DRs_iterative_scheme_z_t}
\end{equation}
and
\begin{equation*}
        z_{t+1} = \argmin_z\,\nabla g_{i_t}(z_t)^T (z-z_t) + \frac{1}{2\lambda}\norm{z-(2x_{t+1}-u_t)}^2,
\end{equation*}
which are called inexact oDRs (ioDRs) and inexact sDRs (isDRs) respectively.

ioDRs or isDRs can also be derived from another point of view, which is based on \emph{proximal gradient}~\cite{nesterov2007gradient}. Here we use ioDRs as an example.
The proximal gradient method uses the proximal mapping of the
nonsmooth part to minimize composite functions~(\ref{eq:composite-form-online})~\cite{lee2013proximal}:
\begin{align*}
   x_{t+1} &= \prox_{\lambda_th}\left(x_t-\lambda_t\nabla g(x_t)\right)
%\\        &= \argmin_y\,t_kh(y) + \frac{1}{2}\norm{y-x_k+t_k\nabla g(x_k)}^2
 \\     & = \argmin_y\,\nabla g(x_t)^T(y-x_t) + \frac{1}{2\lambda_t}\norm{y-x_t}^2 + h(y),
\end{align*}
where $\lambda_t$ denotes the $t$-th step length. %and $G_{t_kf}(x_k)$ is a
%\emph{composite gradient step}.
Then the online PG (OPG) is straight forward, that is at around $t$ solving the following
optimization problem with the linearization of only $t$-th loss function $g_t(x)$~\cite{xiao2010dual,duchi2012dual,suzuki2013dual}:
\begin{align*}
   x_{t+1}   = \argmin_y\,\nabla g_t(x_t)^T(y-x_t) + \frac{1}{2\lambda_t}\norm{y-x_t}^2 + h(y).
\end{align*}
%where $\tilde{g}_t(x) = \frac{1}{t}\sum_{i=1}^t g_i (x)$.
Then the ioDRs can be seen as a combination of OPG with DRs.

\section{Convergence Analysis for DRs}
The procedure of batch DRs is summarized in Algotihtm~\ref{alg:prox-gradient-with-DRs}.
It is clear that $x^*$ is a solution of Problem (\ref{eq:composite-form-online}) if and only if $0\in \partial g(x^*)+\partial h(x^*)$, which is equivalent to $x^* - \lambda\partial g(x^*)\in x^*+\lambda\partial h(x^*)$, where $\lambda>0$. It is clear that $\partial g(x)$ and $\partial h(x)$ are two monotone set-valued operators~\cite{rockafellar1997convex}, and the resolvent operators $R^\lambda_{\partial g} := (I+\lambda\partial g)^{-1}$ and $R^\lambda_{\partial h} :=(I+\lambda\partial h)^{-1}$ are both single valued. Thus if the loss functions $g_t$ are smooth, we have $x^*=R^\lambda_{\partial h}(x^* - \lambda\nabla g(x^*))$, which immediately gives an accuracy measure proposed in~\cite{he2011convergence} of a vector $x$ to a solution of Problem (\ref{eq:composite-form-online}) by
\begin{equation*}
        \varepsilon_g(x,\lambda)=x-R^\lambda_{\partial h}(x - \lambda\nabla g(x)).
\end{equation*}

After $t$ iterations of~(\ref{eq:DRs_iterative_scheme}), from $\varepsilon=x_t-R^\lambda_{\partial h}(x_t - \lambda\nabla g(x_t))$, we have $\frac{\varepsilon}{\lambda}\in \nabla g(x_t)+\partial h(x_t-\varepsilon)$.
Since $g(x)$ and $f(x)$ are convex functions, using their (sub)gradients, we have
\begin{align}
 g(x_t)-g(x^*)&\leq \langle x_t-x^*,\nabla g(x_t)\rangle
 \label{eq:grad_ineq_1}
\\ h(x_t-\varepsilon)-h(x^*)& \leq \langle x_t-\varepsilon-x^*,ph(x_t-\varepsilon)\rangle,
\label{eq:grad_ineq_2}
\\ h(x_t)-h(x_t-\varepsilon)& \leq \langle  \varepsilon,ph(x_t) \rangle ,
\label{eq:grad_ineq_3}
  \end{align}
where $ph(x_t-\varepsilon) \in \partial h(x_t-\varepsilon)$ is the subgradient of $h(x)$ at $x_t-\varepsilon$ satisfying $\frac{\varepsilon}{\lambda} = \nabla g(x_t)+ph(x_t-\varepsilon)$, and $ph(x_t) \in \partial h(x_t)$ is any subgradient of $h(x)$ at $x_t$.

Adding~(\ref{eq:grad_ineq_1}), ~(\ref{eq:grad_ineq_2}), and ~(\ref{eq:grad_ineq_3}) together yields
\begin{align}
 & g(x_t)+h(x_t)-(g(x^*)+h(x^*)) \nonumber \\
 \leq &\langle x_t-x^*,\nabla g(x_t)+ph(x_t-\varepsilon)\rangle+\langle \varepsilon,ph(x_t)-ph(x_t-\varepsilon)\rangle, \label{eq:grad_ineq_4}\\
 \leq &\langle x_t-x^*,\frac{\varepsilon}{\lambda} \rangle+\langle \varepsilon,ph(x_t)-ph(x_t-\varepsilon)\rangle. \nonumber
  \end{align}

If $h(x)$ is Lipschitz continuous with Lipschitz constant $L_h$, then we have $\|ph(x)\|\leq  L_h$, where $ph(x)$ is any subgradient of $h(x)$ at $x$. Furthermore $x_t$ is bounded for its convergence. Thus we have $\|ph(x_t)-ph(x_t-\varepsilon)\|=O(1)$, $\|x_t-x^*\|=O(1)$, thus the following convergence result holds.

\begin{theorem}
  \label{thm:batch_regret_bound}
Assume $g(x)$ is differentiable, $h(x)$ is Lipschitz continuous. Let the sequence $\{x_t,z_t,u_t\}$ be generated by DRs. Then we have
  \begin{align}
 g(x_T)+h(x_T)-(g(x^*)+h(x^*))=O(\varepsilon_g(x_{T},\lambda)). \nonumber
  \end{align}
\end{theorem}

\begin{remark} From  above theorem, we notice that if we have a faster convergent rate of $\varepsilon_g(x_{T},\lambda)$, then we have a faster convergent rate of the optimizing value. It is showed in Theorem 3.1 of~\cite{he2011convergence} that after $t$ iterations of~(\ref{eq:DRs_iterative_scheme}), we have
\begin{align*}
\|\varepsilon_g(x_{t},\lambda)\|^2 = O(1/t).
 \end{align*}
 Thus we have
\end{remark}
\begin{corollary}
 \label{thm:batch_regret_bound}
Assume $g(x)$ is differentiable, $h(x)$ is Lipschitz continuous. Let the sequence $\{x_t,z_t,u_t\}$ be generated by DRs in Algotihtm~\ref{alg:prox-gradient-with-DRs}. Then we have
  \begin{align}
 g(x_T)+h(x_T)-(g(x^*)+h(x^*))=O(1/\sqrt{T}). \nonumber
  \end{align}
\end{corollary}

\begin{remark}
 If $g(x)$ is Lipschitz continuous with Lipschitz constant $L_g$ and the set of all $\partial h(x)$ is bounded by some constant $N_{\partial h}$, then we can obtain an explicit bound by using Lemma 3.1 of~\cite{he2011convergence}, which shows us that $\|(x_t-x^*)+\lambda(\nabla g(x_t)-\nabla g(x^*))\|$ is monotonically decreasing. Thus we have $\|x_t-x^*\| \leq \|x_0-x^*\|+4 \lambda L_g$ from the following inference
  \begin{align}
& \|x_t-x^*\|-2\lambda L_g \nonumber\\
\leq &\|(x_t-x^*)+\lambda(\nabla g(x_t)-\nabla g(x^*))\| \nonumber \\
\leq &\|(x_0-x^*)+\lambda(\nabla g(x_0)-\nabla g(x^*))\|  \label{eq:grad_ineq_5} \\ \leq & \|x_0-x^*\|+2\lambda L_g. \nonumber
  \end{align}
  Further from Theorem 3.1 in~\cite{he2011convergence}, we have
  \begin{align*}
 &\|\varepsilon_g(x_{t},\lambda)\|  \\
\leq &\frac{1}{\sqrt{t+1}} \|(x_t-x^*)+\lambda(\nabla g(x_t)-\nabla g(x^*))\| \\
\leq &\frac{1}{\sqrt{t+1}} (\|x_t-x^*\|+2\lambda L_g)\\
\leq &\frac{1}{\sqrt{t+1}} (\|x_0-x^*\|+6\lambda L_g).
  \end{align*}
Then according  to~(\ref{eq:grad_ineq_4}), we have
  \begin{align}
  &g(x_T)+h(x_T)-(g(x^*)+h(x^*))\nonumber \\
 \leq &\langle x_T-x^*,\frac{\varepsilon}{\lambda} \rangle+\langle \varepsilon,ph(x_T)-ph(x_T-\varepsilon)\rangle \nonumber \\
 \leq &\|x_T-x^*\|\|\frac{\varepsilon}{\lambda}\|+\|\varepsilon\|\|ph(x_T)-ph(x_T-\varepsilon)\rangle\|  \label{eq:grad_ineq_7} \\
  \leq & \frac{1}{\lambda}(\|x_0-x^*\|+4 \lambda L_g)\| \varepsilon\|+2N_{\partial h} \| \varepsilon\|. \nonumber
  \end{align}
 Thus we have the following corollary.
\end{remark}

\begin{corollary}
Assume $g(x)$ is differentiable, both $h(x)$ and $g(x)$ are Lipschitz continuous, and the set of all $\partial h(x)$ is bounded by some constant $N_{\partial h}$. Let the sequence $\{x_t,z_t,u_t\}$ be generated by DRs. Then we have
  \begin{align}
 &g(x_T)+h(x_T)-(g(x^*)+h(x^*)) \nonumber\\
 \leq &\frac{1}{\lambda\sqrt{T+1}}(\|x_0-x^*\|+4 \lambda L_g)( \|x_0-x^*\|+6\lambda L_g)+\frac{2N_{\partial h}}{\sqrt{T+1}} (\|x_0-x^*\|+6\lambda L_g). \nonumber
  \end{align}
\end{corollary}

\begin{remark}
If furthermore $\nabla g(x)$ is Lipschitz continuous with $L_{\nabla g}$, then from~(\ref{eq:grad_ineq_5}) we have $\|x_t-x^*\| \leq \frac{1+\lambda L_{\nabla g}}{1-\lambda L_{\nabla g}}\|x_0-x^*\| $ from the following derivation
 \begin{align}
&(1-\lambda L_{\nabla g})\|x_t-x^*\| \nonumber\\
\leq &\|(x_t-x^*)+\lambda(\nabla g(x_t)-\nabla g(x^*))\| \nonumber\\
\leq &\|(x_0-x^*)+\lambda(\nabla g(x_0)-\nabla g(x^*))\|  \label{eq:grad_ineq_8}
\\ \leq & (1+\lambda L_{\nabla g})\|x_0-x^*\|. \nonumber
  \end{align}
 Then due to the similar formulation, we have $\|\varepsilon_g(x_{t},\lambda)\| \leq \frac{1+\lambda L_{\nabla g}}{\sqrt{t+1}}\|x_0-x^*\|$.
 Put these new inequalities into~(\ref{eq:grad_ineq_7}), then we have the following rate.
\end{remark}

\begin{corollary}
Assume $g(x)$ is differentiable, both $h(x)$ and $g(x)$ are Lipschitz continuous, the set of all $\partial h(x)$ is bounded by some constant $N_{\partial h}$, and $\nabla g(x)$ is Lipschitz continuous with $L_{\nabla g}$. Let the sequence $\{x_t,z_t,u_t\}$ be generated by DRs. Then we have
  \begin{align*}
 &g(x_T)+h(x_T)-(g(x^*)+h(x^*)) \\
 \leq &(\frac{(1+\lambda L_{\nabla g})}{\lambda(1-\lambda L_{\nabla g})}+2N_{\partial h})\frac{1+\lambda L_{\nabla g}}{\sqrt{T+1}}\|x_0-x^*\|.
  \end{align*}
\end{corollary}

These above results are all based and derived from the formulation of~\cite{he2011convergence}, we wonder what we can obtain from the DR iteration formulations~(\ref{eq:DRs_iterative_scheme}),~(\ref{eq:DRs_iterative_scheme_z_t}), and~(\ref{eq:DRs_iterative_scheme_x_t}).
From~(\ref{eq:DRs_iterative_scheme}), we have
  \begin{align}
& 0\in \partial h(x_{t+1})+\frac{1}{\lambda}(x_{t+1}-u_t) \nonumber \\
\Rightarrow &-\frac{1}{\lambda}(x_{t+1}-u_t) \in \partial h(x_{t+1}).\label{eq:grad_ineq_9}
  \end{align}
From~(\ref{eq:DRs_iterative_scheme_z_t}), we have
 \begin{align*}
&0\in \partial g(z_{t+1})+\frac{1}{\lambda}(z_{t+1}-(2x_{t+1}-u_t)) \\
\Rightarrow &-\frac{1}{\lambda}(z_{t+1}-(2x_{t+1}-u_t)) \in \partial g(z_{t+1}).
  \end{align*}
  Using~(\ref{eq:DRs_iterative_scheme_x_t}) yields
 \begin{align}
-\frac{1}{\lambda}(u_{t+1}-x_{t+1}) \in \partial g(z_{t+1}).
 \label{eq:grad_ineq_11}
  \end{align}
Since $h$ and $g$ are convex functions and their subgradients are given in~(\ref{eq:grad_ineq_9}) and~(\ref{eq:grad_ineq_11}) respectively, we have,
\begin{align*}
h(x_{t+1})-h(x^*)\leq &\langle-\frac{1}{\lambda}(x_{t+1}-u_t), x_{t+1}-x^*\rangle,\\
g(z_{t+1})-g(x^*)\leq & \langle-\frac{1}{\lambda}(u_{t+1}-x_{t+1}), z_{t+1}-x^*\rangle
\end{align*}
Adding above together yields
\begin{align}
&h(x_{t+1})+g(z_{t+1})-(h(x^*)+g(x^*))\nonumber \\
\leq &\langle-\frac{1}{\lambda}(x_{t+1}-u_t), x_{t+1}-x^*\rangle+ \langle-\frac{1}{\lambda}(u_{t+1}-x_{t+1}), z_{t+1}-x^*\rangle \label{eq:grad_ineq_13}\\
= &\frac{1}{\lambda}[u_{t+1}(x^*-z_{t+1})+x_{t+1}(z_{t+1}-x_{t+1})+u_t(x_{t+1}-x^*)] \nonumber
\end{align}
If $g(x)$ is Lipschitz continuous with Lipschitz constant $L_g$, then we have $\|g(x_{t+1})-g(z_{t+1})\|\leq L_g\|x_{t+1}-z_{t+1}\|$, adding with~(\ref{eq:grad_ineq_13}) yields
\begin{align*}
&h(x_{t+1})+g(x_{t+1})-(h(x^*)+g(x^*)) \\
\leq &\frac{1}{\lambda}[u_{t+1}(x^*-z_{t+1})+u_t(x_{t+1}-x^*)]+(L_g+\frac{1}{\lambda}\|x_{t+1}\|)\|x_{t+1}-z_{t+1}\|.
\end{align*}
Thus we have
\begin{theorem}
  \label{thm:batch_regret_bound_2}
Assume $g(x)$ is Lipschitz continuous. Let the sequence $\{x_t,z_t,u_t\}$ be generated by DRs. Then we have
  \begin{align}
 g(x_T)+h(x_T)-(g(x^*)+h(x^*))=O(x_{T}-x^*). \nonumber
  \end{align}
\end{theorem}

\section{Online and stochastic Douglas-Rachford splitting method}
In this section, we generalize the DRs to online and stochastic settings.
The procedure of batch DRs, oDRs, ioDRs, sDRs, and isDRs are summarized in Algotihtm~\ref{alg:prox-gradient-with-DRs},~\ref{alg:online-prox-gradient-with-DRs},
~\ref{alg:inexact_online-prox-gradient-with-DRs},~\ref{alg:sDRs} and~\ref{alg:isDRs} respectively, where $f_1(x)=g_1(x)+h(x)$.

\begin{algorithm}[t]
\caption{A generic DRs}
\label{alg:prox-gradient-with-DRs}
%\begin{AlgorithmSteps}[4]
%Recovering of Low-rank Component from Audio Segments via robust PCA.

\textbf{Input}: starting point $x_0\in\dom (g+h)$.

%\textbf{Initialize}: $D\in \mathbb{R}^{m\times n\times p}, \beta>0, k=0.$

1: \textbf{for} $t=0,1,\cdots,T$ \textbf{do}

2: $x_{t+1} = \argmin_x\,h(x) + \frac{1}{2\lambda}\norm{x-u_t}^2$.

3: $z_{t+1}  = \argmin_z\,g_{t+1}(z) + \frac{1}{2\lambda}\norm{z-(2x_{t+1}-u_t)}^2$.

4: $ u_{t+1}  = u_t + \lambda_t(z_{t+1}-x_{t+1})$.

5: \textbf{end for}

\textbf{Output}: $x_{T+1}$.
\end{algorithm}

\begin{algorithm}[t]
\caption{A generic oDRs}
\label{alg:online-prox-gradient-with-DRs}
%\begin{AlgorithmSteps}[4]
%Recovering of Low-rank Component from Audio Segments via robust PCA.

\textbf{Input}: starting point $x_0\in\dom f_1$.

%\textbf{Initialize}: $D\in \mathbb{R}^{m\times n\times p}, \beta>0, k=0.$

1: \textbf{for} $t=0,\cdots,T$ \textbf{do}

2: $x_{t+1} = \argmin_x\,h(x) + \frac{1}{2\lambda}\norm{x-u_t}^2$.

3: $z_{t+1}  = \argmin_z\,g_{t+1}(z) + \frac{1}{2\lambda}\norm{z-(2x_{t+1}-u_t)}^2$.

4: $ u_{t+1}  = u_t + \lambda_t(z_{t+1}-x_{t+1})$.

5: \textbf{end for}

\textbf{Output}: $x_{T+1}$.
\end{algorithm}

\begin{algorithm}[t]
\caption{A generic ioDRs}
\label{alg:inexact_online-prox-gradient-with-DRs}
%\begin{AlgorithmSteps}[4]
%Recovering of Low-rank Component from Audio Segments via robust PCA.

\textbf{Input}: starting point $x_0\in\dom f_1$.

%\textbf{Initialize}: $D\in \mathbb{R}^{m\times n\times p}, \beta>0, k=0.$

1: \textbf{for} $t=0,\cdots,T$ \textbf{do}

2: $x_{t+1} = \argmin_x\,h(x) + \frac{1}{2\lambda}\norm{x-u_t}^2$.

3: $z_{t+1}  = \argmin_z\,\nabla g_{t+1}(z_t)^T (z-z_t) + \frac{1}{2\lambda}\norm{z-(2x_{t+1}-u_t)}^2$.

4: $ u_{t+1}  = u_t + \lambda_t(z_{t+1}-x_{t+1})$.

5: \textbf{end for}

\textbf{Output}: $x_{T+1}$.
\end{algorithm}

\begin{algorithm}[t]
\caption{A generic sDRs}
\label{alg:sDRs}

\textbf{Input}: starting point $x_0\in\dom f_1$.

%\textbf{Initialize}: $D\in \mathbb{R}^{m\times n\times p}, \beta>0, k=0.$

1: \textbf{for} $t=0,\cdots,T$ \textbf{do}

2: $x_{t+1} = \argmin_x\,h(x) + \frac{1}{2\lambda}\norm{x-u_t}^2$.

3: Randomly select index $i_t$ from the set $\{1,...,T\}$ and solve the subproblem:
\begin{align*}
z_{t+1}  = \argmin_z\,g_{i_{t+1}}(z) + \frac{1}{2\lambda}\norm{z-(2x_{t+1}-u_t)}^2.
\end{align*}

4: $ u_{t+1}  = u_t + \lambda_t(z_{t+1}-x_{t+1})$.

5: \textbf{end for}

\textbf{Output}: $x_{T+1}$.
\end{algorithm}

\begin{algorithm}[t]
\caption{A generic ioDRs}
\label{alg:isDRs}
%\begin{AlgorithmSteps}[4]
%Recovering of Low-rank Component from Audio Segments via robust PCA.

\textbf{Input}: starting point $x_0\in\dom f_1$.

%\textbf{Initialize}: $D\in \mathbb{R}^{m\times n\times p}, \beta>0, k=0.$

1: \textbf{for} $t=0,\cdots,T$ \textbf{do}

2: $x_{t+1} = \argmin_x\,h(x) + \frac{1}{2\lambda}\norm{x-u_t}^2$.

3: Randomly select index $i_t$ from the set $\{1,...,T\}$ and solve the subproblem:
\begin{align*}
z_{t+1}  = \argmin_z\,\nabla g_{i_{t+1}}(z_t)^T (z-z_t) + \frac{1}{2\lambda}\norm{z-(2x_{t+1}-u_t)}^2.
\end{align*}

4: $ u_{t+1}  = u_t + \lambda_t(z_{t+1}-x_{t+1})$.

5: \textbf{end for}

\textbf{Output}: $x_{T+1}$.
\end{algorithm}

\subsection{Regret Analysis for oDRs}

The goal of oDRs is to achieve low regret w.r.t. a static predictor on a sequence of functions
\begin{equation*}
f_T(x) = \frac{1}{T}\sum_{t=1}^T g_t (x) + h(x).
\end{equation*}
According to Algotihtm~\ref{alg:online-prox-gradient-with-DRs}, formally, at every round of the algorithm we make a prediction $x_t$ and then receive the function $f_t(x)=g_t (x) + h(x)$. That is at round $t-1$, we obtain $x_{t}$ by solving the following problem:
\begin{equation*}
        x_t^*=\argmin_{x}\,g_t (x) + h(x)
\end{equation*}
 with only single DR iteration based on the warm start $x_{t-1}$. In batch optimization we set $f_t=f$ for all $t$ while in stochastic optimization we choose $f_t$ to be the average of some random subset of $\{f_1,...,f_T\}$.

In this work, we seek bounds on the standard regret in the online learning setting with respect to $x^*$, defined as
\begin{equation*}
        R(T,x^*):=\frac{1}{T}\sum_{t=1}^T (g_t (x_t) + h(x_t))- [\frac{1}{T}\sum_{t=1}^T g_t (x^*) + h(x^*)]
\end{equation*}

As pointed by~(\ref{eq:grad_ineq_8}) and Theorem 3.1 in~\cite{he2011convergence}, with the notation $\varepsilon_{g_t}(x_t,\lambda)=x_t-R^\lambda_{\partial h}(x_t- \lambda\nabla g_t(x_t))$ in mind, we have in each iteration that
 \begin{align*}
\|x_t-x_t^*\| \leq \frac{1+\lambda L_{\nabla g_t}}{1-\lambda L_{\nabla g_t}}\|x_{t-1}-x_t^*\|
\end{align*}
 and
 \begin{align*}
\|\varepsilon_{g_t}(x_{t},\lambda)\|^2\leq \frac{1}{2}\|(x_{t-1}-x_t^*)+\lambda(\nabla g_t(x_{t-1})-\nabla g_t(x_t^*))\|^2, \nonumber
\end{align*}
which means that
 \begin{align*}
\frac{O(1)}{\lambda}\in \nabla g_t(x_t)+\partial h(x_t-O(1)).
\end{align*}
Following the same procedure as~(\ref{eq:grad_ineq_4}), we have
\begin{align*}
 &g_t(x_t)+h(x_t)-(g_t(x^*)+h(x^*)) \\
 \leq &\langle x_t-x^*,\frac{O(1)}{\lambda} \rangle+\langle O(1),ph(x_t)-ph(x_t-O(1))\rangle =O(1).
  \end{align*}
  Summing up above formulas for $t\in \{1,...,T\}$, we obtain the following result:
\begin{theorem}
  \label{thm:online_regret_bound}
Assume all $g_t(x)$ are differentiable, $h(x)$ and $g_t(x)$ are Lipschitz continuous, the set of all $\partial h(x)$ is bounded by some constant $N_{\partial h}$, $\nabla g_t(x)$ is Lipschitz continuous with $L_{\nabla g_t}$, all $L_{\nabla g_t}$ and $x_t^*$ are bounded. Let the sequence $\{x_t,z_t,u_t\}$ be generated by oDRs. Then we have
  \begin{align}
 R(T,x^*)=O(1). \nonumber
  \end{align}
\end{theorem}

\subsection{Convergence analysis of sDRs}

From~(\ref{eq:grad_ineq_8}) and Theorem 3.1 in~\cite{he2011convergence}, we have in each iteration of~Algotihtm~\ref{alg:sDRs} that
 \begin{align*}
\|x_t-x_t^*\| \leq \frac{1+\lambda L_{\nabla g_{i_t}}}{1-\lambda L_{\nabla g_{i_t}}}\|x_{t-1}-x_t^*\|
\end{align*}
 and
 \begin{align*}
\|\varepsilon_{g_{i_t}}(x_{t},\lambda)\|^2\leq \frac{1}{2}\|(x_{t-1}-x_t^*)+\lambda(\nabla g_{i_t}(x_{t-1})-\nabla g_{i_t}(x_t^*))\|^2, \nonumber
\end{align*}
which means that
 \begin{align*}
\frac{O(1)}{\lambda}\in \nabla g_{i_t}(x_t)+\partial h(x_t-O(1)).
\end{align*}
Following the same procedure as~(\ref{eq:grad_ineq_4}), we have
\begin{align*}
 &g_{i_t}(x_t)+h(x_t)-(g_{i_t}(x^*)+h(x^*)) \\
 \leq &\langle x_t-x^*,\frac{O(1)}{\lambda} \rangle+\langle O(1),ph(x_t)-ph(x_t-O(1))\rangle =O(1).
  \end{align*}
Thus we obtain the following result:
\begin{theorem}
  \label{thm:stochastic_convergence_rate}
Assume all $g_t(x)$ are differentiable, $h(x)$ and $g_t(x)$ are Lipschitz continuous, the set of all $\partial h(x)$ is bounded by some constant $N_{\partial h}$, $\nabla g_t(x)$ is Lipschitz continuous with $L_{\nabla g_t}$, all $L_{\nabla g_t}$ and $x_t^*$ are bounded. Let the sequence $\{x_t,z_t,u_t\}$ be generated by~Algotihtm~\ref{alg:sDRs}. Then we have
  \begin{align}
 g(x_T)+h(x_T)-(g(x^*)+h(x^*))=O(1/\sqrt{T}). \nonumber
  \end{align}
\end{theorem}

\section{Computational experiments}
In this section, we demonstrate the performance of oDRs and sDRs in solving several machine learning problems. We present simulation results to show the convergence of the objective in oDRs and sDRs. We also compare them with batch DRs and OADM~\cite{wang2012online}. We set $\lambda_t=1$ for all the updates of $u_{t+1}$, and $\lambda=[0.1, 1, 10, 20]$. All the experiments show that oDRs and sDRs outperform OADM.
\subsection{Lasso}
The lasso problem is formulated as follows:
\begin{align}
  \minimize_{x \in \reals^{n\times 1}}\,\frac{1}{T}\sum_{t=1}^T \|a_t^Tx-b_t\|^2 + \mu\norm{x}_1,
  \label{eq:lasso}
\end{align}
where $a_t, x\in \reals^{n\times 1}$ and $b_t$ is a scalar. The three updatas of DRs are:
\begin{align*} %\label{eq:DRs_iterative_scheme_lasso}
  x_{t+1}
        = &\argmin_x\,\mu\norm{x}_1 + \frac{1}{2\lambda}\norm{x-u_t}^2 \\
        =&\sign(u_t)\cdot\max\{\abs{u_t}-\mu\lambda,0\},
        \\ z_{t+1}  = &\argmin_z\,\frac{1}{T}\sum_{t=1}^T \|a_t^Tz-b_t\|^2 + \frac{1}{2\lambda}\norm{z-(2x_{t+1}-u_t)}^2
        \\   = &(A^TA+\frac{T}{2\lambda}I)^{-1}[A^Tb+\frac{T}{2\lambda}(2x_{t+1}-u_t)],
        %\label{eq:DRs_iterative_scheme_z_t}
 \\   u_{t+1}   =& u_t + \lambda_t(z_{t+1}-x_{t+1}),
          %   \label{eq:DRs_iterative_scheme_x_t}
\end{align*}
where $A=(a_1, ..., a_T)$ and $b = (b_1, ..., b_T)^T$.
The differences of oDRs and ioDRs from DRs is the update of $z_{t+1}$, which are:
\begin{align*}
        z_{t+1}& = \argmin_z\,\|a_t^Tz-b_t\|^2 + \frac{1}{2\lambda}\norm{z-(2x_{t+1}-t_t)}^2  \\
        &= [a_t^Ta_t+\frac{1}{2\lambda}I]^{-1}[a_tb_t+\frac{1}{2\lambda}(2x_{t+1}-u_t)]
\end{align*}
and
\begin{align*}
        z_{t+1}& = \argmin_z\,2(a_t^Tz_t-b_t)a_t^T (z-z_t) + \frac{1}{2\lambda}\norm{z-(2x_{t+1}-u_t)}^2 \\
        & = 2x_{t+1}-u_t - 2\lambda(a_t^Tz_t-b_t)a_t \nonumber
\end{align*}
respectively.

Our experiments mainly follow the lasso example in~\cite{wang2012online}. We first randomly generated $A$ with 1000 examples of
dimensionality 100. $A$ is then normalized along the columns. Then, a true $x_0$ is randomly generated
with certain sparsity pattern for lasso, and we set the number of nonzeros as 10. $b$ is calculated
by adding Gaussian noise to $Ax_0/T$, where $T=1000$ is number of examples. We set $\mu=0.1\times \|A^Tb/T\|_\infty$ and $\eta=1$ in OADM~\cite{wang2012online}. All experiments are implemented in Matlab.

\begin{figure}[!h]
\subfigure[$\lambda$=0.1]{
\includegraphics[width=0.5\textwidth]{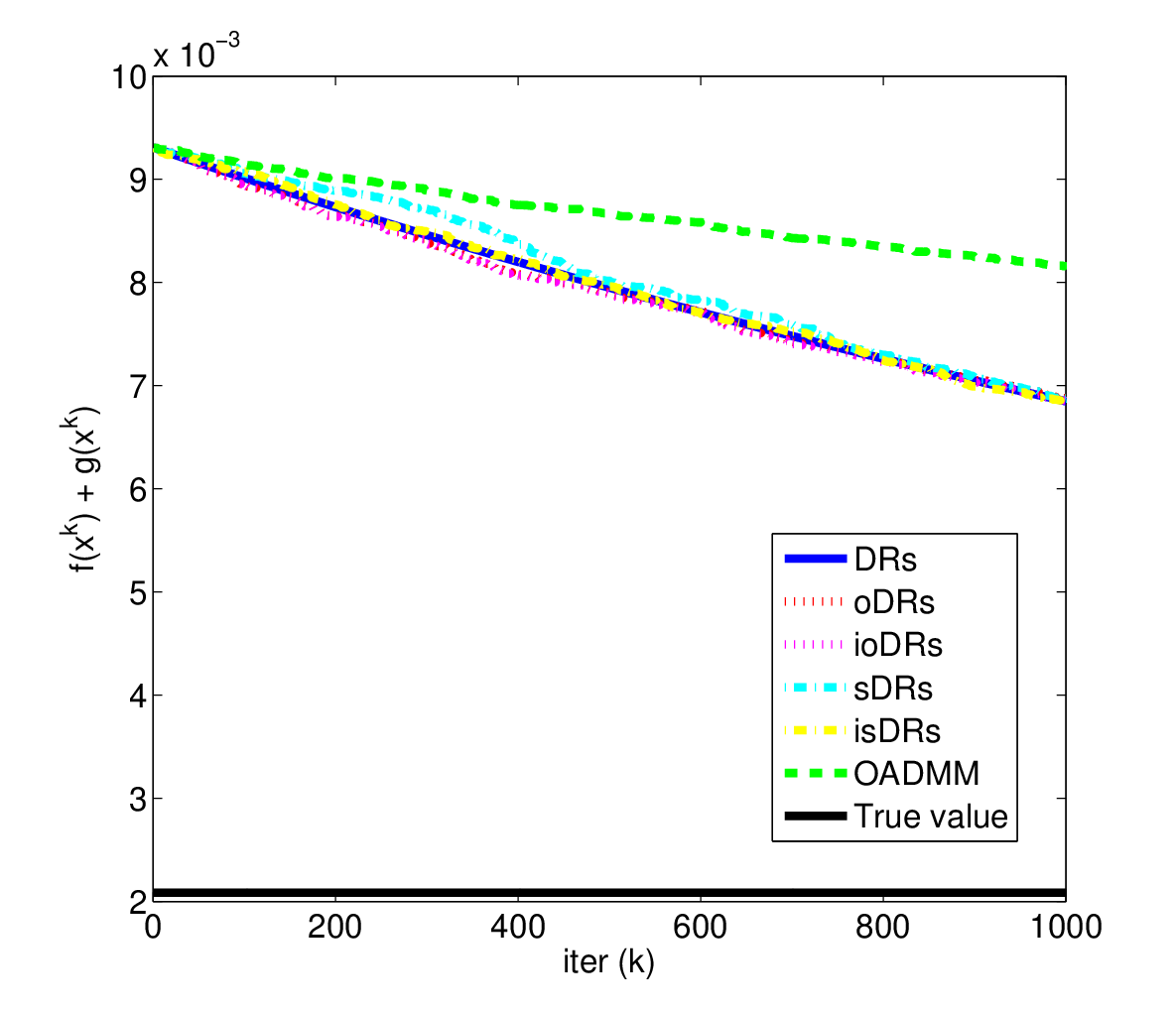}
}
\subfigure[$\lambda$=1]{
\includegraphics[width=0.5\textwidth]{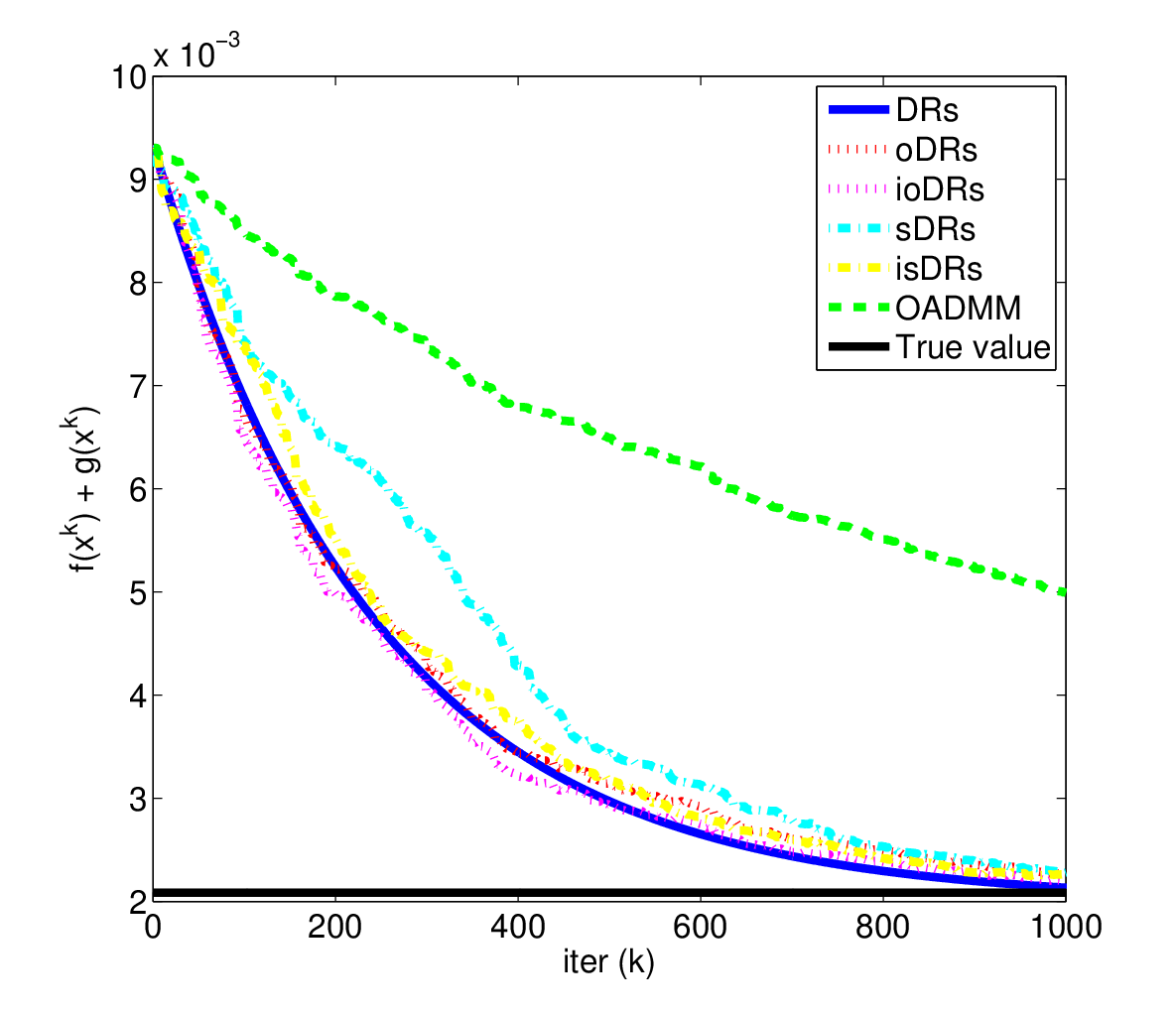}
}
\subfigure[$\lambda$=10]{
\includegraphics[width=0.5\textwidth]{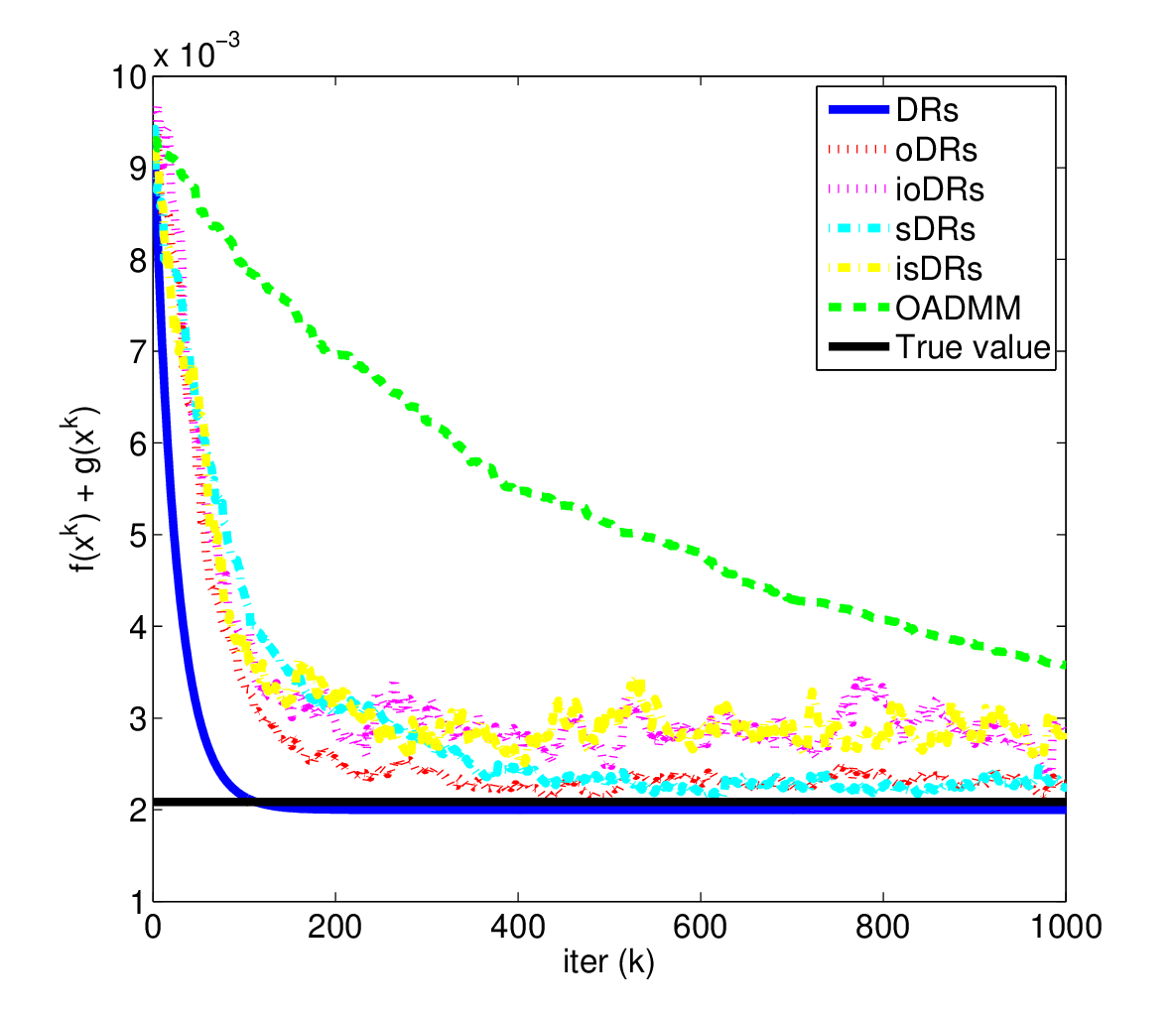}
}
\subfigure[$\lambda$=20]{
\includegraphics[width=0.5\textwidth]{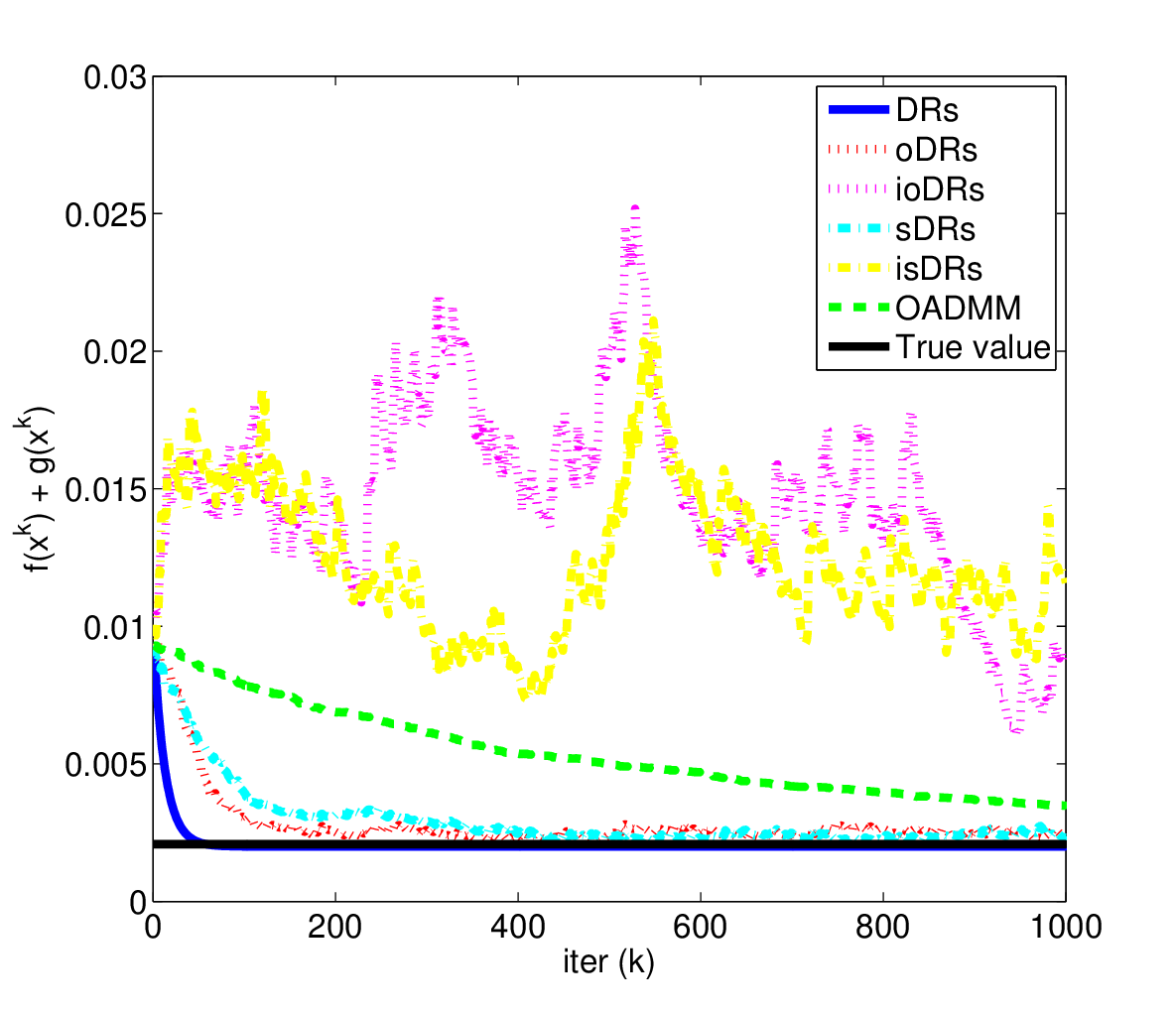}
}
\caption{ The convergence of objective value in DRs, oDRs, ioDRs, sDRs, isDRs, OADM, and real objective value for the lasso problem~(\ref{eq:lasso}). }
\label{fig:Lasso_Objecti_Value}       % Give a unique label
\end{figure}

In Figure~\ref{fig:Lasso_Objecti_Value}, the objective value of the problem is depicted against the iteration times. In this example, oDRs and sDRs show faster convergence than OADM, ioDRs, and isDRs. The main reason for the slow convergence of ioDRs or isDRs is because the linearization of the objective in each iteration~(\ref{eq:inexact_Online_DRs_iterative_scheme_z_t}). We observe that OADM takes even longer iterations to achieve a certain precision, although the regret bound is more tighter than the bound obtained in this work for oDRs ($O(1/\sqrt{T})$ in OADM~\cite{wang2012online}, while $O(1)$ in oDRs). Thus we believe and conjecture that the regret $R(T,x^*)$ in Theorem~\ref{thm:online_regret_bound} is indeed $O(1/\sqrt{T})$.

\subsection{Logistic regression}

The logistic regression problem is formulated as follows:
\begin{align}
  \minimize_{w \in \reals^n}\,\frac{1}{T}\sum_{i=1}^T \log (1+\exp(-y_i w^Tx_i)) + \mu\norm{w}_1,
  \label{eq:logi}
\end{align}
where $x^{(1)},\dots,x^{(T)}$ are samples with labels
$y^{(1)},\dots,y^{(T)}\in\{0,1\}$, the regularization term $\norm{w}_1$ promotes sparse solutions
and $\mu$ balances goodness-of-fit and sparsity.

The three updatas of DRs are:
\begin{align*} %\label{eq:DRs_iterative_scheme_lasso}
  w_{t+1}
        &= \argmin_w\,\mu\norm{w}_1 + \frac{1}{2\lambda}\norm{w-u_t}^2 \\
        &=\sign(w_t)\cdot\max\{\abs{w_t}-\mu\lambda,0\},
        \\ z_{t+1} & = \argmin_z\,\frac{1}{T}\sum_{i=1}^T \log (1+\exp(-y_i z^Tx_i)) + \frac{1}{2\lambda}\norm{z-(2w_{t+1}-u_t)}^2,
        %\label{eq:DRs_iterative_scheme_z_t}
 \\   u_{t+1}  & = u_t + \lambda_t(z_{t+1}-x_{t+1}).
          %   \label{eq:DRs_iterative_scheme_x_t}
\end{align*}

The differences of oDRs and ioDRs from DRs is the update of $z_{t+1}$, which are:
\begin{align*} %\label{eq:DRs_iterative_scheme_lasso}
z_{t+1} & = \argmin_z\, \log (1+\exp(-y_t z^Tx_t)) + \frac{1}{2\lambda}\norm{z-(2w_{t+1}-u_t)}^2,
\end{align*}
and
\begin{align*} %\label{eq:DRs_iterative_scheme_lasso}
      z_{t+1} & = \argmin_z\, -\frac{1+\exp(-y_t z_t^Tx_t)}{\exp(-y_t z_t^Tx_t)}y_tx_t(z-z_t) + \frac{1}{2\lambda}\norm{z-(2w_{t+1}-u_t)}^2.
\end{align*}
The sDRs and isDRs are almost the same, we will not repeat here.

\subsection{Total variation minimization}

The total variation (TV) minimization problem is formulated as follows:
\begin{align}
  \minimize_{x \in \reals^{n\times 1}}\,\frac{1}{T}\sum_{t=1}^T \|a_t^Tx-b_t\|^2 + \mu\sum_{t=1}^{T-1}\norm{x_{t+1}-x_t}_1,
  \label{eq:tv}
\end{align}
where $a_t, x\in \reals^{n\times 1}$ and $b_t$ is a scalar.

The three updatas of DRs are:
\begin{align*} %\label{eq:DRs_iterative_scheme_lasso}
  x_{t+1}
        &= \argmin_x\,\mu\norm{Dx}_1 + \frac{1}{2\lambda}\norm{x-u_t}^2 \\
        & = \argmin_x\,\mu\norm{Dx}_1 + \frac{1}{2\lambda}\norm{EDx-EDu_t}^2,
        \\ z_{t+1} & = \argmin_z\,\frac{1}{T}\sum_{t=1}^T \|a_t^Tz-b_t\|^2 + \frac{1}{2\lambda}\norm{z-(2x_{t+1}-u_t)}^2
        \\  & = (A^TA+\frac{T}{2\lambda}I)^{-1}[A^Tb+\frac{T}{2\lambda}(2x_{t+1}-u_t)],
        %\label{eq:DRs_iterative_scheme_z_t}
 \\   u_{t+1}  & = u_t + \lambda_t(z_{t+1}-x_{t+1}),
          %   \label{eq:DRs_iterative_scheme_x_t}
\end{align*}
where $A=(a_1, ..., a_T)$, $b = (b_1, ..., b_T)^T$ and $D$ is an upper bi-diagonal matrix with diagonal 1 and off-diagonal -1, and $ED=I$. Here we should notice that the update of $x$ is obtained by solving a small scale lasso problem. Here we employ the DR algorithm to solve it.
The differences of oDRs and ioDRs from DRs is the update of $z_{t+1}$, which are:
\begin{align*}
        z_{t+1}& = \argmin_z\,\|a_t^Tz-b_t\|^2 + \frac{1}{2\lambda}\norm{z-(2x_{t+1}-t_t)}^2 \\
         & = [a_t^Ta_t+\frac{1}{2\lambda}I]^{-1}[a_tb_t+\frac{1}{2\lambda}(2x_{t+1}-u_t)]
\end{align*}
and
\begin{align*}
        z_{t+1}& = \argmin_z\,2(a_t^Tz_t-b_t)a_t^T (z-z_t) + \frac{1}{2\lambda}\norm{z-(2x_{t+1}-u_t)}^2 \\
        & = 2x_{t+1}-u_t - 2\lambda(a_t^Tz_t-b_t)a_t
\end{align*}
respectively.

Based on the numerical tests, we propose the following two very credible conjectures:

\begin{conjecture}
Assume $g(x)$ is differentiable, $h(x)$ is Lipschitz continuous. Let the sequence $\{x_t,z_t,u_t\}$ be generated by ioDRs. Then we have
  \begin{align}
 R(T,x^*)=O(1). \nonumber
  \end{align}
\end{conjecture}

and

\begin{conjecture}
Assume $g(x)$ is differentiable, $h(x)$ is Lipschitz continuous. Let the sequence $\{x_t,z_t,u_t\}$ be generated by isDRs. Then we have
  \begin{align}
 g(x_T)+h(x_T)-(g(x^*)+h(x^*))=O(1/\sqrt{T}). \nonumber
  \end{align}
\end{conjecture}

\section{Conclusion}

In this paper, we propose the efficient online and stochastic learning algorithm
named online DRs (oDRs) and stochastic DRs (sDRs) respectively. New proof techniques
have been developed to analyze the regret of DRs, oDRs and convergence of sDRs,
which shows that DRs and oDRs have $O(1/\sqrt{T})$ and $O(1)$ regret respectively, and sDRs enjoys a convergence rate of $O(1/\sqrt{T})$. Finally, we illustrate
the efficiency of oDRs and sDRs in solving several machine learning problems.

\frenchspacing
\bibliographystyle{plain}
\bibliography{OnlineDRs_arxiv}

\begin{thebibliography}{10}

\bibitem{beck2003mirror}
Amir Beck and Marc Teboulle.
\newblock Mirror descent and nonlinear projected subgradient methods for convex
  optimization.
\newblock {\em Operations Research Letters}, 31(3):167--175, 2003.

\bibitem{combettes2009iterative}
Patrick~L Combettes.
\newblock Iterative construction of the resolvent of a sum of maximal monotone
  operators.
\newblock {\em J. Convex Anal}, 16(4):727--748, 2009.

\bibitem{combettes2011proximal}
Patrick~L Combettes and Jean-Christophe Pesquet.
\newblock Proximal splitting methods in signal processing.
\newblock In {\em Fixed-Point Algorithms for Inverse Problems in Science and
  Engineering}, pages 185--212. Springer, 2011.

\bibitem{douglas1956numerical}
Jim Douglas and HH~Rachford.
\newblock On the numerical solution of heat conduction problems in two and
  three space variables.
\newblock {\em Transactions of the American mathematical Society},
  82(2):421--439, 1956.

\bibitem{duchi2010composite}
John Duchi, Shai Shalev-Shwartz, Yoram Singer, and Ambuj Tewari.
\newblock Composite objective mirror descent.
\newblock 2010.

\bibitem{duchi2012dual}
John~C Duchi, Alekh Agarwal, and Martin~J Wainwright.
\newblock Dual averaging for distributed optimization.
\newblock In {\em Communication, Control, and Computing (Allerton), 2012 50th
  Annual Allerton Conference on}, pages 1564--1565. IEEE, 2012.

\bibitem{eckstein1992douglas}
Jonathan Eckstein and Dimitri~P Bertsekas.
\newblock On the douglas¡ªrachford splitting method and the proximal point
  algorithm for maximal monotone operators.
\newblock {\em Mathematical Programming}, 55(1-3):293--318, 1992.

\bibitem{gabay1976dual}
Daniel Gabay and Bertrand Mercier.
\newblock A dual algorithm for the solution of nonlinear variational problems
  via finite element approximation.
\newblock {\em Computers \& Mathematics with Applications}, 2(1):17--40, 1976.

\bibitem{he2011convergence}
BS~He and XM~Yuan.
\newblock On convergence rate of the douglas-rachford operator splitting
  method.
\newblock {\em Mathematical Programming, under revision}, 2011.

\bibitem{lee2013proximal}
Jason~D Lee, Yuekai Sun, and Michael~A Saunders.
\newblock Proximal newton-type methods for minimizing composite functions.
\newblock 2013.

\bibitem{lions1979splitting}
Pierre-Louis Lions and Bertrand Mercier.
\newblock Splitting algorithms for the sum of two nonlinear operators.
\newblock {\em SIAM Journal on Numerical Analysis}, 16(6):964--979, 1979.

\bibitem{nesterov2007gradient}
Yurii Nesterov.
\newblock Gradient methods for minimizing composite objective function, 2007.

\bibitem{nesterov2009primal}
Yurii Nesterov.
\newblock Primal-dual subgradient methods for convex problems.
\newblock {\em Mathematical programming}, 120(1):221--259, 2009.

\bibitem{rockafellar1997convex}
R~Tyrell Rockafellar.
\newblock {\em Convex analysis}, volume~28.
\newblock Princeton university press, 1997.

\bibitem{shi2013audio}
Ziqiang Shi, Jiqing Han, Tieran Zheng, and Shiwen Deng.
\newblock Audio segment classification using online learning based tensor
  representation feature discrimination.
\newblock {\em IEEE transactions on audio, speech, and language processing},
  21(1-2):186--196, 2013.

\bibitem{suzuki2013dual}
Taiji Suzuki.
\newblock Dual averaging and proximal gradient descent for online alternating
  direction multiplier method.
\newblock In {\em Proceedings of the 30th International Conference on Machine
  Learning (ICML-13)}, pages 392--400, 2013.

\bibitem{tibshirani1996regression}
Robert Tibshirani.
\newblock Regression shrinkage and selection via the lasso.
\newblock {\em Journal of the Royal Statistical Society. Series B
  (Methodological)}, pages 267--288, 1996.

\bibitem{wang2012online}
Huahua Wang and Arindam Banerjee.
\newblock Online alternating direction method.
\newblock {\em arXiv preprint arXiv:1206.6448}, 2012.

\bibitem{xiao2010dual}
Lin Xiao.
\newblock Dual averaging methods for regularized stochastic learning and online
  optimization.
\newblock {\em The Journal of Machine Learning Research}, 11:2543--2596, 2010.

\bibitem{zinkevich2003online}
Martin Zinkevich.
\newblock Online convex programming and generalized infinitesimal gradient
  ascent.
\newblock 2003.

\end{thebibliography}

\end{document}